\documentclass[11pt,paper]{amsart}

\usepackage{amssymb,latexsym}
\usepackage[matrix,arrow]{xy}
\usepackage[usenames]{color}
\usepackage{slashed}
\usepackage{parskip}
 \numberwithin{equation}{section}
\newtheorem{theorem}{Theorem}[section]
\newtheorem{corollary}{Corollary}[section]
\newtheorem{lemma}{Lemma}[section]
\newtheorem{proposition}{Proposition}[section]

\theoremstyle{definition}

\theoremstyle{remark}
\newtheorem{remark}{Remark}[section]

\def\proof{{\indent\bf Proof.}}
\def\endproof{\hfill\hbox{$\sqcup$}\llap{\hbox{$\sqcap$}}\medskip}

\begin{document}
\title[Eigenvalues of the Laplacian]{Eigenvalues of the Laplacian
on \\
 Riemannian manifolds* }
\author [Q. -M. Cheng and X. Qi]{ Qing-Ming Cheng  and Xuerong Qi}
\address{Qing-Ming Cheng\\  Department of Mathematics, Graduate School of Science and Engineering,
Saga University, Saga 840-8502,  Japan, cheng@cc.saga-u.ac.jp}
\address{Xuerong Qi\\  Department of Mathematics,  Graduate School  of Science and Engineering,
Saga University, Saga 840-8502,  Japan, qixuerong609@gmail.com}
\renewcommand{\thefootnote}{\fnsymbol{footnote}}


\footnotetext{* Research partially supported by a Grant-in-Aid for
Scientific Research from JSPS.}


\begin{abstract}

For a bounded domain $\Omega$ with a piecewise smooth boundary  in a
complete Riemannian manifold $M$, we study eigenvalues of the
Dirichlet eigenvalue problem of the Laplacian. By making use of a
fact that eigenfunctions  form  an orthonormal  basis of
$L^2(\Omega)$ in place of the Rayleigh-Ritz formula, we obtain
inequalities for eigenvalues of the Laplacian. In particular, for
lower order eigenvalues, our results extend the results of Chen and
Cheng \cite{CC}.
\end{abstract}
\maketitle
\renewcommand{\sectionmark}[1]{}

\section{introduction}

Let $\Omega \subset M$ be a bounded domain with a piecewise smooth boundary
 $\partial \Omega$    in an
 $n$-dimensional complete Riemannian manifold $M$. We consider the
 following Dirichlet eigenvalue problem of the Laplacian:
\begin{equation}
\left \{ \aligned \Delta u&=-\lambda u \quad \text{in \ \ $\Omega$}, \\
u&=0 \quad \ \ \ \ \text{on $\partial \Omega$}. \endaligned \right.
\end{equation}
 It is well known that the spectrum of this problem is   real
 and  discrete:
 $$
 0<\lambda_1< \lambda_2\leq \lambda_3\leq \cdots\nearrow\infty,
$$
where each $\lambda_i$ has finite multiplicity which is repeated
according to its multiplicity.

When $M$ is an $n$-dimensional Euclidean space
$\mathbf{R}^n$,
   Payne, P\'olya and Weinberger \cite{PPW1} proved
\begin{equation}
\lambda_{k+1} - \lambda_{k} \leq \frac{4}{kn} \sum_{i=1}^{k} \lambda_{i} .
\end{equation}
Hile and Protter \cite{HP} generalized the above result to
\begin{equation}
\sum_{i=1}^{k} \frac{\lambda_{i}}{\lambda_{k+1} - \lambda_{i}} \geq \frac{kn}{4}.
\end{equation}
In 1991, a much  sharper  inequality was obtained by Yang \cite{Y} (cf. \cite{CY4}):
\begin{equation}
\sum_{i=1}^{k}(\lambda_{k+1} - \lambda_{i} )^{2}
\leq \frac{4}{n} \sum_{i=1}^{k}(\lambda_{k+1} - \lambda_{i} ) \lambda_{i}.
\end{equation}
 When $M$ is an $n$-dimensional unit sphere $S^n(1)$,
Cheng and Yang \cite{CY1} have proved an optimal universal
inequality:
\begin{equation}
\sum_{i=1}^{k}(\lambda_{k+1} - \lambda_{i} )^{2} \leq \frac{4}{n}
 \sum_{i=1}^{k}(\lambda_{k+1} - \lambda_{i} ) (\lambda_{i}+\displaystyle{\frac{n^2}{4}}).
\end{equation}

For the Dirichlet eigenvalue problem of the Laplacian on a bounded domain in
an $n$-dimensional complete Riemannian manifold $M$, Chen and Cheng
\cite{CC} and El Soufi, Harrell and Ilias \cite{EHI} have proved,
independently,
\begin{equation}
\sum_{i=1}^{k}(\lambda_{k+1} - \lambda_{i} )^{2} \leq \frac{4}{n}
 \sum_{i=1}^{k}(\lambda_{k+1} - \lambda_{i} )
 (\lambda_{i}+\displaystyle{\frac{n^2}{4}H_0^2}),
\end{equation}
where $H_0^2$ is a nonnegative constant which only depends on $M$
and $\Omega$. When $M$ is the unit sphere,  $H_0^2=1$, the above
inequality is best possible, which becomes the result of Cheng and
Yang  \cite{CY1}. For the Dirichlet eigenvalue problem of the
Laplacian on a bounded domain in a hyperbolic space, universal
inequalities for eigenvalues have been obtained by Cheng and Yang
\cite{CY5}. For complex projective spaces and so on, see \cite{CCC},
\cite{CY1}  and \cite{CY3}.

For lower order eigenvalues of the eigenvalue problem (1.1),  when $M$ is the Euclidean space $\mathbf R^n$,
the following conjecture of  Payne, P\'olya and Weinberger is well
known:

\noindent {\bf Conjecture of PPW.}  \ \ {\it For a bounded domain
$\Omega$ in  $\mathbf R^n$,   eigenvalues of the eigenvalue problem
{\rm (1.1)} satisfy
\begin{enumerate}
\item  $\dfrac{\lambda_2}{\lambda_1}
\leq \left.\dfrac{\lambda_2}{\lambda_1}\right|_{\mathbf{B}^{n}}
 =\dfrac{j^2_{n/2,1}}{j^2_{n/2-1,1}}$,
\item  $\dfrac{\lambda_2+\lambda_3+\cdots
+\lambda_{n+1}}{\lambda_1}\leq
n\dfrac{j^2_{n/2,1}}{j^2_{n/2-1,1}}$,
\end{enumerate}
\noindent where $\mathbf{B}^{n}$ is the $n$-dimensional unit ball in
$\mathbf{R}^{n}$, $j_{p,k}$ denotes the $k$-th positive zero of the
standard Bessel function $J_{p}(x)$ of the first kind of order $p$.}

 For the conjecture (1) of Payne, P\'olya and Weinberger, many
mathematicians studied it. For examples, Payne, P\'olya and
Weinberger \cite {PPW1}, Brands \cite{B}, de Vries \cite{dV},  Chiti
\cite{Chi}, Hile and Protter \cite{HP}, Marcellini \cite{Ma}  and so
on. Finally, Ashbaugh and Benguria \cite{AB2} (cf. \cite{AB1} and
\cite{AB3}) solved this conjecture.

For the conjecture (2) of Payne, P\'olya and Weinberger, when $n=2$,
Brands \cite{B} improved the bound
$\dfrac{\lambda_2+\lambda_3}{\lambda_1}\leq 6$ of Payne, P\'olya and
Weinberger \cite{PPW1}, he proved
$\dfrac{\lambda_2+\lambda_3}{\lambda_1}\leq 3+\sqrt7$. Furthermore,
Hile and Protter \cite{HP}  obtained
$\dfrac{\lambda_2+\lambda_3}{\lambda_1}\leq 5.622$. In \cite{Ma},
Marcellini proved $\dfrac{\lambda_2+\lambda_3}{\lambda_1}\leq
(15+\sqrt{345})/6$. Recently, Chen and Zheng \cite{CZ} have proved
$\dfrac{\lambda_2+\lambda_3}{\lambda_1}\leq 5.3507$.
 For a
general dimension $n\geq 2$, Ashbaugh and Benguria \cite{AB4} proved
 \begin{equation}
 \frac{\lambda_2+\lambda_3+\cdots +\lambda_{n+1}}{\lambda_1}\leq
 n +4.
 \end{equation}
 Furthermore,  Ashbaugh and Benguria \cite{AB4} (cf. Hile and Protter \cite{HP} )
improved the above result to
 \begin{equation}
 \frac{\lambda_2+\lambda_3+\cdots +\lambda_{n+1}}{\lambda_1}\leq
 n+3 +\dfrac{\lambda_1}{\lambda_2}.
 \end{equation}
Very recently, Cheng and Qi \cite{CQ} have proved that, for any
$1\leq j\leq n+2$, eigenvalues satisfy at least one of the
following:
$$
\dfrac{\lambda_2}{\lambda_1}<2-\dfrac{\lambda_1}{\lambda_j},
$$
$$\frac{\lambda_2+\lambda_3+\cdots +\lambda_{n+1}}{\lambda_1}\leq
 n+3 +\dfrac{\lambda_1}{\lambda_j}.
 $$

When $M$ is the  $n$-dimensional unit sphere $S^{n}(1)$, that is,
for a  bounded domain $\Omega$ in  $S^{n}(1)$, Cheng, Sun and Yang
\cite {SCY} have proved
\begin{equation}
\dfrac{\lambda_2+\lambda_3+\cdots + \lambda_{n+1}}{\lambda_1}
\leq
n+4+\dfrac{n^2}{\lambda_1}.
\end{equation}
For a  general complete Riemannian manifold $M$, Chen and Cheng
\cite{CC} have proved that there exists a non-negative constant
$H_0$ such that
\begin{equation}
\dfrac{\lambda_2+\lambda_3+\cdots + \lambda_{n+1}}{\lambda_1}\leq
n+4+\dfrac{n^2H_0^2}{\lambda_1}.
\end{equation}
In this paper, by making use of the fact that eigenfunctions form an
orthonormal  basis of $L^2(\Omega)$ in place of the Rayleigh-Ritz
formula, we obtain inequalities for eigenvalues of the Laplacian. In
particular, we  improve the above result.


\section{Estimates for lower order eigenvalues}

In this section, first of all,
we will mainly focus our mind on the investigation for lower order
eigenvalues of the Dirichlet eigenvalue problem of the Laplacian by
making use of the fact that eigenfunctions form an orthonormal
basis of $L^2(\Omega)$ in place of the Rayleigh-Ritz formula. We
prove the following:

\begin{theorem} Let $M$ be an $n$-dimensional complete Riemannian manifold,
$\Omega \subset M $  a bounded domain with a piecewise smooth boundary
 $\partial \Omega$. Then,  the lower order eigenvalues of
the Dirichlet eigenvalue problem of the Laplacian
satisfy
\begin{equation*}
\begin{aligned}
&\dfrac{\lambda_2+\lambda_3+\cdots + \lambda_{n+1}}{\lambda_1}\\
\leq \ & n+\sqrt{\biggl(\dfrac{n^2H_0^2}{\lambda_1}+4\biggl)
\dfrac{(2-\dfrac{\lambda_1}{\lambda_2})\dfrac{n^2H_0^2}{\lambda_1}+3+\dfrac{\lambda_1}
{\lambda_2}+\sqrt{\big(3+\dfrac{\lambda_1}{\lambda_2}
+\dfrac{n^2H_0^2}{\lambda_2}\big)^2+4(1-\dfrac{\lambda_1}{\lambda_2})\dfrac{n^2H_0^2}{\lambda_2}}}{2}}\\
\end{aligned}
\end{equation*}
where $H_0$ is a non-negative constant depending on $M$ and $\Omega$ only.
\end{theorem}

\begin{remark}
It is not hard to prove, from $ \dfrac{\lambda_1}{\lambda_2}<1$,
$$
\dfrac{(2-\dfrac{\lambda_1}{\lambda_2})\dfrac{n^2H_0^2}{\lambda_1}
+3+\dfrac{\lambda_1}{\lambda_2}+\sqrt{\big(3+\dfrac{\lambda_1}{\lambda_2}
+\dfrac{n^2H_0^2}{\lambda_2}\big)^2+4(1-\dfrac{\lambda_1}{\lambda_2})\dfrac{n^2H_0^2}{\lambda_2}}}{2}<\dfrac{n^2H_0^2}{\lambda_1}+4.
$$
\end{remark}

In particular, when $M$ is an $n$-dimensional complete minimal submanifold in the Euclidean space $\mathbf R^N$,
we have
\begin{corollary}
Let $\Omega$ be a  bounded domain in an $n$-dimensional complete
minimal submanifold  $M$  in $\mathbf R^N$. Then,  we have
\begin{equation*}
\dfrac{\lambda_2+\lambda_3+\cdots + \lambda_{n+1}}{\lambda_1} \leq n+2\sqrt{3+\dfrac{\lambda_1}{\lambda_2}}.
\end{equation*}
\end{corollary}

Since $M$ is a complete Riemannian manifold, from a
theorem of Nash, there exists an isometric immersion $\varphi: M\to
\mathbf R^N$ from
 $M$ into  a Euclidean space $\mathbf R^N$.
 Let $(x^1,\cdots,x^n)$ denote an arbitrary local coordinate system of  $M$.
 For any point $p\in \Omega$,  we can write  $\varphi(p)=(y_1, y_2, \cdots, y_N)$
 with
$$
y_{\alpha}=y_{\alpha}(x^1,\cdots,x^n), \qquad 1\leq \alpha \leq N,
$$
which is the position vector of $p $ in  $\mathbf R^N$. Thus, we have
 $$
 g_{ij}=g\Big(\frac{\partial}{\partial x^i},\frac{\partial}{\partial x^j}\Big)
 =\Big\langle \sum_{\alpha=1}^N\frac{\partial y_\alpha}{\partial x^i}\frac{\partial}{\partial y_\alpha},
 \sum_{\beta=1}^N\frac{\partial y_\beta}{\partial x^j}\frac{\partial }{\partial y_\beta} \Big\rangle
 =\sum_{\alpha=1}^N\frac{\partial y_{\alpha}}{\partial x^i}
 \frac{\partial y_{\alpha}}{\partial x^j},
 $$
 where $g$ denotes the induced metric of $M$ from $\mathbf R^N$,
 $\langle \ , \ \rangle$ is the standard inner product  in $\mathbf R^N$.
We denote  the gradient of  a function $f$  by $\nabla f$. Then, the
following lemma holds, which is proved by  Chen and Cheng \cite{CC}.

 \begin{lemma}
 $$\aligned
  \sum_{\alpha=1}^N g(\nabla y_{\alpha}, \nabla y_\alpha)&=\sum_{\alpha=1}^N |\nabla
  y_{\alpha}|^{2}=n,\\
\sum_{\alpha=1}^{N}(\Delta y_\alpha)^2&=n^2 |H|^2,\\
\sum_{\alpha=1}^{N}\Delta y_\alpha \nabla y_{\alpha} &=0,
\endaligned$$ and for any function $u\in C^\infty(M)$,
 $$ \ \ \ \ \ \ \ \ \ \
 \sum_{\alpha=1}^{N}\Big(g(\nabla y_{\alpha},
\nabla u)\Big)^2 =\sum_{\alpha=1}^{N}\Big(\nabla y_{\alpha}\cdot
\nabla u\Big)^2=|\nabla u|^2,
$$
where $|H|$ is the mean curvature of  $M$.
\end{lemma}

\vskip 2mm
\noindent
{\it Proof of Theorem 2.1.}
Let  $u_j$ be the  eigenfunction corresponding
to the eigenvalue $\lambda_j$ such that $\{u_j\}_{j=1}^{\infty}$ becomes
an orthonormal basis of $L^2(\Omega)$. Hence,  $\int_{\Omega} u_i u_j=\delta_{ij} $ for $\forall\
i, j=1, 2, \cdots$.
Defining
$$
a_{\alpha j}=\int_{\Omega} y_{\alpha} u_1u_{j+1},
$$
since $u_1$ does not change sign in $\Omega$, we can assume $u_1>0$ in $\Omega$.
We consider  the $N\times N$-matrix $A=(a_{\alpha j})$.
From the orthogonalization of Gram and Schmidt, there exist an upper triangle matrix $R=(R_{\alpha j})$ and  an orthogonal matrix $Q=(q_{\alpha \beta})$ such that
$R=QA$. Thus,
$$
R_{\alpha j}=\sum_{\beta=1}^Nq_{\alpha\beta}a_{\beta j}
=\int_{\Omega} \sum_{\beta=1}^Nq_{\alpha\beta}y_{\beta} u_1u_{j+1}=0, \ \text{for} \ 1\leq j <\alpha\leq N.
$$
Defining $\bar y_{\alpha}=\sum_{\gamma=1}^Nq_{\alpha\gamma}y_{\gamma}$, we have
$$
\int_{\Omega} \bar y_{\alpha} u_1u_{j+1}=\int_{\Omega} \sum_{\gamma=1}^Nq_{\alpha\gamma}y_{\gamma} u_1u_{j+1}=0,
\ \text{for} \ 1\leq j <\alpha\leq N.
$$
Putting
$$
z_{\alpha}=\bar y_{\alpha}-b_{\alpha}, \  b_{\alpha}=\int_{\Omega} \bar y_{\alpha}u_1^2,
\  \text{for}   \ \ 1\leq \alpha \leq N
$$
and
$$
A_{\alpha j}=\int_{\Omega}z_{\alpha}u_1u_{j},
$$
we have
\begin{equation}
 \ \ \ \ \ \ \ \ A_{\alpha j}=0,\qquad \mbox{for}\ 1\leq j\leq \alpha \leq N.
\end{equation}
Defining
$$
B_{\alpha j}=\int_{\Omega}u_j \nabla z_\alpha \cdot \nabla u_1 \ \ \
\ \
$$
and
$$
C_{\alpha j}=\int_{\Omega}u_ju_1\Delta z_\alpha, \ \ \ \ \ \ \ \ \
$$
from the Stokes theorem, we obtain
\begin{equation*}
\begin{aligned}
-\lambda_jA_{\alpha j}&=\int_{\Omega}z_{\alpha}u_1\Delta u_j=\int_{\Omega}\Delta(z_\alpha u_1)u_j\\
&=\int_{\Omega}\biggl( 2\nabla z_\alpha \cdot\nabla u_1-
\lambda_1z_\alpha u_1+u_1\Delta z_\alpha \biggl)u_j\\
&=-\lambda_1A_{\alpha j}+2B_{\alpha j}+C_{\alpha j},
\end{aligned}
\end{equation*}
namely,
\begin{equation}
2B_{\alpha j}=(\lambda_1-\lambda_j)A_{\alpha j}-C_{\alpha j}.
\end{equation}

Since $\{u_j\}_{j=1}^{\infty}$ is an orthonormal basis in
$L^2(\Omega)$ and $A_{\alpha j}=0,\ {\rm for}\  1\leq j\leq \alpha
\leq N$,  we have
\begin{equation}
z_{\alpha}u_1=\sum_{j=\alpha+1}^{\infty}A_{\alpha j}u_j \ \ {\rm
and} \quad
  \|z_{\alpha}u_1\|^2=\sum_{j=\alpha+1}^{\infty}A_{\alpha j}^2.
\end{equation}
Furthermore,
\begin{equation}
\int_{\Omega}u_1^2z_{\alpha}\Delta z_{\alpha}=\sum_{j=\alpha+1}^{\infty}A_{\alpha j}C_{\alpha j},
\end{equation}
\begin{equation}
2\int_{\Omega}z_{\alpha}u_1 \nabla z_{\alpha}\cdot \nabla u_1
=2\sum_{j=\alpha+1}^{\infty}A_{\alpha j}B_{\alpha j}
=\sum_{j=\alpha+1}^{\infty}(\lambda_1-\lambda_j)A_{\alpha
j}^2-\sum_{j=\alpha+1}^{\infty}A_{\alpha j}C_{\alpha j}.
\end{equation}
Since for any function $f\in C^2(\Omega)\cap C(\bar\Omega)$,
\begin{equation}
-2\int_{\Omega}fu_1 \nabla  f \cdot \nabla u_1
=\int_{\Omega}u_1^2f\Delta f+\int_{\Omega}|\nabla  f|^2u_1^2,
\end{equation}
 we have
\begin{equation}
\int_{\Omega}|\nabla  z_{\alpha}|^2u_1^2=
-\int_{\Omega}z_{\alpha}u_1\biggl(2\nabla  z_{\alpha}\cdot \nabla
u_1+u_1\Delta z_{\alpha}\biggl).
\end{equation}
We obtain
\begin{equation}
\begin{aligned}
\sum_{j=\alpha+1}^{\infty}(\lambda_j-\lambda_1)A_{\alpha
j}^2=\int_{\Omega}|\nabla  z_{\alpha}|^2u_1^2.
\end{aligned}
\end{equation}
For any positive integer $k$,  we have
\begin{equation*}
\begin{aligned}
\sum_{j=\alpha+1}^{\infty}(\lambda_j-\lambda_1)A_{\alpha j}^2
=&\sum_{j=\alpha+1}^{k}(\lambda_j-\lambda_1)A_{\alpha j}^2+\sum_{j=k+1}^{\infty}(\lambda_j-\lambda_1)A_{\alpha j}^2\\
\geq &\sum_{j=\alpha+1}^{k}(\lambda_j-\lambda_1)A_{\alpha j}^2
+(\lambda_{k+1}-\lambda_1)\sum_{j=k+1}^{\infty}A_{\alpha j}^2\\
=&\sum_{j=\alpha+1}^{k}(\lambda_j-\lambda_1)A_{\alpha j}^2
+(\lambda_{k+1}-\lambda_1)\sum_{j=\alpha+1}^{\infty}A_{\alpha j}^2-(\lambda_{k+1}-\lambda_1)\sum_{j=\alpha+1}^{k}A_{\alpha j}^2\\
=&\sum_{j=\alpha+1}^{k}(\lambda_j-\lambda_{k+1})A_{\alpha j}^2
+(\lambda_{k+1}-\lambda_1)\sum_{j=\alpha+1}^{\infty}A_{\alpha j}^2.
\end{aligned}
\end{equation*}
Thus, we infer
\begin{equation}
\begin{aligned}
(\lambda_{k+1}-\lambda_1) \|z_{\alpha}u_1\|^2\leq
\sum_{j=\alpha+1}^{k}(\lambda_{k+1}-\lambda_j)A_{\alpha
j}^2+\int_{\Omega}|\nabla  z_{\alpha}|^2u_1^2,
\end{aligned}
\end{equation}
and, in particular,
\begin{equation}
\begin{aligned}
(\lambda_{\alpha+1}-\lambda_1)
\|z_{\alpha}u_1\|^2\leq\int_{\Omega}|\nabla  z_{\alpha}|^2u_1^2.
\end{aligned}
\end{equation}

For any $\alpha$, we have
\begin{equation}
  |\nabla  z_{\alpha}|^2\leq 1.
\end{equation}
 In fact, for any fixed point $p_0\in \Omega$,
 we can choose a new coordinate system $\widetilde y=(\widetilde y_1,\cdots,\widetilde  y_N)$ of $\mathbf R^N$
given by  $\varphi(p)-\varphi(p_0)=\widetilde y(p) B$ such that
 $\frac {\partial }{\partial \widetilde y_1}|_{p_0},\cdots,
 \frac {\partial }{\partial \widetilde y_n}|_{p_0}$ span $T_{p_0}M$
 and at $p_0$, $g\Big(\frac{\partial}{\partial \widetilde y_i},\frac{\partial}{\partial \widetilde y_j}\Big)
 =\delta_{ij}$,
 where  $B=(b_{\alpha \beta})\in O(N)$ is an $N\times N$ orthogonal matrix.
\begin{equation}
\begin{aligned}
  |\nabla  z_{\alpha}|^2(p_0)&=g(\nabla  z_{\alpha}, \nabla  z_{\alpha})\\
&= \sum_{\beta,\gamma=1}^Nq_{\alpha\gamma}q_{\alpha\beta}g(\nabla  y_{\gamma}, \nabla  y_{\beta})\\
&=\sum_{\beta,\gamma=1}^Nq_{\alpha\gamma}q_{\alpha\beta}
g(\sum_{\mu=1}^Nb_{\gamma\mu}\nabla  \widetilde y_{\mu}, \sum_{\nu=1}^Nb_{\beta\nu}\nabla  \widetilde y_{\nu})\\
&=\sum_{\beta,\gamma, \mu, \nu=1}^N
q_{\alpha\gamma}b_{\gamma\mu}q_{\alpha\beta}b_{\beta\nu}g(\nabla  \widetilde y_{\mu}, \nabla  \widetilde y_{\nu})\\
&=\sum_{j=1}^n\bigl(\sum_{\beta=1}^Nq_{\alpha\beta}b_{\beta j}\bigl)^2\leq 1,
\end{aligned}
\end{equation}
since $QB$ is an orthogonal matrix when $B$ and $Q$ are orthogonal
matrices. Therefore, (2.11) holds because $p_0$ is an arbitrary
point. Since Lemma 2.1 also holds for $z_{\alpha}$ from the
definition of $z_{\alpha}$, for any positive constant $t>\frac12$,
we have, from Lemma 2.1 and (2.11),
\begin{equation}
\begin{aligned}
  &\sum_{\alpha=1}^N(\lambda_{\alpha+1}-\lambda_1)\int_{\Omega}|\nabla  z_{\alpha}|^2u_1^{t+1}\\
  \geq &\sum_{j=1}^n(\lambda_{j+1}-\lambda_1)\int_{\Omega}|\nabla  z_j|^2u_1^{t+1}+
 (\lambda_{n+1}-\lambda_1) \sum_{A=n+1}^N\int_{\Omega}|\nabla  z_A|^2u_1^{t+1}\\
 =& \sum_{j=1}^n(\lambda_{j+1}-\lambda_1)\int_{\Omega}|\nabla  z_j|^2u_1^{t+1}+
 (\lambda_{n+1}-\lambda_1) \int_{\Omega}(n-\sum_{j=1}^n|\nabla  z_j|^2)u_1^{t+1}\\
 =& \sum_{j=1}^n(\lambda_{j+1}-\lambda_1)\int_{\Omega}|\nabla  z_j|^2u_1^{t+1}+
 (\lambda_{n+1}-\lambda_1) \int_{\Omega}\sum_{j=1}^n(1-|\nabla  z_j|^2)u_1^{t+1}\\
\geq &\sum_{j=1}^n(\lambda_{j+1}-\lambda_1)\int_{\Omega}|\nabla
z_j|^2u_1^{t+1}+
  \int_{\Omega}\sum_{j=1}^n(\lambda_{j+1}-\lambda_1)(1-|\nabla  z_j|^2)u_1^{t+1}\\
=& \sum_{j=1}^n(\lambda_{j+1}-\lambda_1)\int_{\Omega}u_1^{t+1}.
\end{aligned}
\end{equation}
On the other hand, from the  Stokes theorem  and  the
Cauchy-Schwarz inequality, we obtain
\begin{equation}
\begin{aligned}
\int_{\Omega}|\nabla  z_{\alpha}|^2u_1^{t+1} =&-\int_{\Omega}
z_{\alpha}u_1\biggl (u_1^t \Delta z_{\alpha}
+(1+t)u_1^{t-1}\nabla  z_{\alpha}\cdot \nabla u_1\biggl)\\
\leq &~\|z_{\alpha}u_1\|\cdot \|u_1^t\Delta
z_{\alpha}+(1+t)u_1^{t-1}\nabla z_{\alpha}\cdot \nabla u_1\|,
\end{aligned}
\end{equation}
and
\begin{equation}
\begin{aligned}
\int_{\Omega}|\nabla  z_{\alpha}|^2u_1^{2} =&-\int_{\Omega}
z_{\alpha}u_1 \biggl (u_1\Delta z_{\alpha}
+2\nabla  z_{\alpha}\cdot  \nabla u_1\biggl)\\
\leq &~\|z_{\alpha}u_1\| \cdot\|u_1\Delta z_{\alpha}+2\nabla
z_{\alpha}\cdot \nabla u_1\|.
\end{aligned}
\end{equation}
From (2.10),  (2.13), (2.14) and (2.15), we derive
\begin{equation}
\begin{aligned}
&\sum_{j=1}^n(\lambda_{j+1}-\lambda_1)\int_{\Omega}u_1^{t+1}\\
\leq &
\sum_{\alpha=1}^N(\lambda_{\alpha+1}-\lambda_1)\int_{\Omega}|\nabla
z_{\alpha}|^2u_1^{t+1}\\ \leq &
\sum_{\alpha=1}^N\dfrac{\int_{\Omega}|\nabla
z_{\alpha}|^2u_1^{2}}{\|z_{\alpha}u_1\| ^2} \int_{\Omega}|\nabla
z_{\alpha}|^2u_1^{t+1}\\
 \leq &\sum_{\alpha=1}^N\|u_1\Delta
z_{\alpha}+2\nabla  z_{\alpha}\cdot \nabla u_1\|\cdot \|u_1^t\Delta
z_{\alpha}+(1+t)u_1^{t-1}\nabla  z_{\alpha}\cdot \nabla u_1\|\\
 \leq &\sqrt{\sum_{\alpha=1}^N\|u_1\Delta z_{\alpha}+2\nabla
z_{\alpha}\cdot \nabla u_1\|^2 \cdot\sum_{\alpha=1}^N\|u_1^t\Delta
z_{\alpha}+(1+t)u_1^{t-1}\nabla  z_{\alpha}\cdot \nabla u_1\|^2}~.
\end{aligned}
\end{equation}

Since Lemma 2.1 also holds for $z_{\alpha}$ from the
definition of $z_{\alpha}$, we have
\begin{equation}
\begin{aligned}
\sum_{\alpha=1}^N\|u_1\Delta z_{\alpha}+2\nabla
z_{\alpha}\cdot\nabla u_1\|^2 =&\int_{\Omega}(n^2|H|^2u_1^2+4|\nabla
u_1|^2)\\ \leq &~ n^2\sup_{\Omega} |H|^2+4\lambda_1
\end{aligned}
\end{equation}
and
\begin{equation}
\begin{aligned}
\sum_{\alpha=1}^N\|u_1^t\Delta z_{\alpha}+(1+t)u_1^{t-1}\nabla
z_{\alpha}\cdot \nabla u_1\|^2
=&\int_{\Omega}\biggl(n^2|H|^2u_1^{2t}+\dfrac{(1+t)^2}{t^2}|\nabla
u_1^t|^2\biggl)\\ \leq & \biggl(n^2\sup_{\Omega}
|H|^2+\dfrac{(1+t)^2}{2t-1}\lambda_1\biggl)\int_{\Omega}u_1^{2t}.
\end{aligned}
\end{equation}
Putting (2.17) and (2.18) into (2.16), we obtain
\begin{equation}
\begin{aligned}
&\sum_{j=1}^n(\lambda_{j+1}-\lambda_1)\leq
B(t)\sqrt{\biggl(n^2\sup_{\Omega}
|H|^2+4\lambda_1\biggl)\biggl(n^2\sup_{\Omega}
|H|^2+\dfrac{(1+t)^2}{2t-1}\lambda_1\biggl)} \ ,
\end{aligned}
\end{equation}
where
$$
B(t)=\dfrac{\sqrt{\int_{\Omega}u_1^{2t}}}{\int_{\Omega}u_1^{t+1}}.
$$
Since the spectrum of the Dirichlet eigenvalue problem of the
Laplacian is an  invariant of  isometries, we know that the above
inequality holds for any isometric immersion from $M$ into a
Euclidean space. Now we define $\Phi$ by
$$\Phi=\big\{\varphi; \varphi  \text{\ is an isometric immersion  from $M$ into a Euclidean space}\big\}.
$$
Defining
$$
H_0^2=\inf_{\varphi\in\Phi}\sup_{\Omega}|H|^2,
$$
we have
\begin{equation}
\begin{aligned}
&\sum_{j=1}^n(\lambda_{j+1}-\lambda_1)\leq
B(t)\sqrt{(n^2H_0^2+4\lambda_1)\biggl(n^2H_0^2+\dfrac{(1+t)^2}{2t-1}\lambda_1\biggl)}~.
\end{aligned}
\end{equation}

Next, we need to estimate $B(t)$ as a function of $t$
by making use of the same method as Brands \cite{B}.
Let $u=u_1^t-u_1\int_{\Omega}u_1^{t+1}$. We know that $u$ is a trial function
for $\lambda_2$. Hence, we have
$$
\lambda_2\leq \dfrac{\int_{\Omega}|\nabla u|^2}{\int_{\Omega}u^2}.
$$
According to a direct calculation, we obtain
$$
\dfrac{\lambda_2}{\lambda_1}\leq \dfrac{\dfrac{t^2}{2t-1}B(t)^2-1}{B(t)^2-1}
$$
since
$B(t)^2-1=\dfrac{\int_{\Omega}u^{2}}{\left(\int_{\Omega}u_1^{t+1}\right)^{2}}>0$
for $t>1$. Let $a=\dfrac{\lambda_2}{\lambda_1}>1$. We have
$$
\left(a-\dfrac{t^2}{2t-1}\right)B(t)^2\leq a-1.
$$
When $1<t<a+\sqrt{a^2-a}$, we can infer
$$
B(t)\leq \sqrt{\dfrac{(a-1)(2t-1)}{a(2t-1)-t^2}}~.
$$
Therefore, we obtain
\begin{equation}
\begin{aligned}
&\sum_{j=1}^n(\lambda_{j+1}-\lambda_1)\leq
\sqrt{\left(n^2H_0^2+4\lambda_1\right)\left(n^2H_0^2+\dfrac{(1+t)^2}{2t-1}\lambda_1\right)\dfrac{(a-1)(2t-1)}{a(2t-1)-t^2}}~.
\end{aligned}
\end{equation}

Letting $b=\dfrac{n^2H_0^2}{\lambda_1}$ and defining  a  function
\begin{equation}
f(t)=\dfrac{b(2t-1)+(1+t)^2}{a(2t-1)-t^2},
\end{equation}
we have
\begin{equation}
\begin{aligned}
&\sum_{j=1}^n(\lambda_{j+1}-\lambda_1)\leq
\sqrt{\lambda_1(a-1)\left(n^2H_0^2+4\lambda_1\right)f(t)}.
\end{aligned}
\end{equation}
If we take $t=1$, $f(1)=\dfrac{b+4}{a-1}$. Thus, we obtain  the result of Chen and Cheng \cite{CC}.
Furthermore, we try  to get the minimum of $f(t)$ under $1\leq t\leq a+\sqrt{a^2-a}$.
It is not difficult to prove  that  the minimum of $f(t)$
is attained at
$$
t_0=\dfrac{a+b-1+\sqrt{(a+b-1)^2+8a(a+b+1)}}{2(a+b+1)}.
$$
Since $g(s)=t_0$ as a function of $s=a+b$, is a decreasing function of $s$ in the
interval $[a, \infty)$, we have
$$
1=g(\infty)<t_0\leq g(a)<a+\sqrt{a^2-a}.
$$
By a direct computation, we have
\begin{equation*}
\begin{aligned}
a(2t_0-1)-t_0^2 &=\dfrac{2\sqrt{(a+b-1)^2+8a(a+b+1)}
}{4(a+b+1)^2}\\
&\times \biggl(2a(a+b+1)-(a+b-1)-\sqrt{(a+b-1)^2+8a(a+b+1)}\biggl).
\end{aligned}
\end{equation*}
 From $\big\{3(a+b)+1\big\}^2-8b(a+b+1)=(a+b-1)^2+8a(a+b+1)$, we get
\begin{equation*}
\begin{aligned}
b(2t_0-1)+(1+t_0)^2 &=\dfrac{2\sqrt{(a+b-1)^2+8a(a+b+1)}
}{4(a+b+1)^2}\\
&\times \biggl(2b(a+b+1)+3(a+b)+1+\sqrt{(a+b-1)^2+8a(a+b+1)}\biggl).
\end{aligned}
\end{equation*}
Thus, we have
\begin{equation}
\begin{aligned}
f(t_0)
&=\dfrac{2b(a+b+1)+3(a+b)+1+\sqrt{(a+b-1)^2+8a(a+b+1)}}{2a(a+b+1)-(a+b-1)-\sqrt{(a+b-1)^2+8a(a+b+1)}}\\
&=\dfrac{2(a+b+1)^2\big\{(2a-1)b+3a+1+\sqrt{(a+b-1)^2+8a(a+b+1)}\big\}}{4a(a+b+1)^2(a-1)}\\
&=\dfrac{(2a-1)b+3a+1+\sqrt{(a+b-1)^2+8a(a+b+1)}}{2a(a-1)}.\\
\end{aligned}
\end{equation}
From (2.23) and (2.24), we obtain
\begin{equation*}
\begin{aligned}
&\sum_{j=1}^n(\lambda_{j+1}-\lambda_1)\\ \leq &~
\sqrt{\lambda_1\left(n^2H_0^2+4\lambda_1\right)\dfrac{(2a-1)b+3a+1+\sqrt{(a+b-1)^2+8a(a+b+1)}}{2a}}\\
=&~\lambda_1\sqrt{\left(\dfrac{n^2H_0^2}{\lambda_1}+4\right)
\dfrac{(2-\dfrac{\lambda_1}{\lambda_2})
\dfrac{n^2H_0^2}{\lambda_1}+3+\dfrac{\lambda_1}{\lambda_2}
+\sqrt{\big(1+\dfrac{n^2H_0^2}{\lambda_2}-\dfrac{\lambda_1}{\lambda_2}\big)^2
+8(1+\dfrac{n^2H_0^2}{\lambda_2}+\dfrac{\lambda_1}{\lambda_2})}}{2}}\\
=&~\lambda_1\sqrt{\left(\dfrac{n^2H_0^2}{\lambda_1}+4\right)\dfrac{(2-\dfrac{\lambda_1}{\lambda_2})\dfrac{n^2H_0^2}{\lambda_1}
+3+\dfrac{\lambda_1}{\lambda_2}+\sqrt{\big(3+\dfrac{\lambda_1}{\lambda_2}
+\dfrac{n^2H_0^2}{\lambda_2}\big)^2+4(1-\dfrac{\lambda_1}{\lambda_2})\dfrac{n^2H_0^2}{\lambda_2}}}{2}}\\
\end{aligned}
\end{equation*}
because $a=\dfrac{\lambda_2}{\lambda_1}$ and
$b=\dfrac{n^2H_0^2}{\lambda_1}$. This finishes the proof of Theorem
2.1. \endproof

 For any positive integer $k$,  we have from (2.9)
\begin{equation}
\begin{aligned}
\frac{\lambda_{k+1}-\lambda_1}{\sum\limits_{j=\alpha+1}^{k}(\lambda_{k+1}-\lambda_j)A_{\alpha
j}^2+\int_{\Omega}|\nabla  z_{\alpha}|^2u_1^2} \leq
\frac{1}{\|z_{\alpha}u_1\|^2}.
\end{aligned}
\end{equation}
From (2.14) and (2.15), we obtain
\begin{equation}
\begin{aligned}
&\dfrac{(\lambda_{k+1}-\lambda_1)\int_{\Omega}|\nabla
z_{\alpha}|^2u_1^{t+1} \int_{\Omega}|\nabla z_{\alpha}|^2u_1^{2}}
{ \sum\limits_{j=\alpha+1}^{k}(\lambda_{k+1}-\lambda_j)A_{\alpha j}^2+\int_{\Omega}|\nabla  z_{\alpha}|^2u_1^2}\\
=&~\dfrac{(\lambda_{k+1}-\lambda_1)\int_{\Omega}|\nabla
z_{\alpha}|^2u_1^{t+1} } {1+
\sum\limits_{j=\alpha+1}^{k}(\lambda_{k+1}-\lambda_j)\dfrac{A_{\alpha
j}^2}
{\int_{\Omega}|\nabla z_{\alpha}|^2u_1^{2}}}\\
\leq &~\|u_1^t\Delta z_{\alpha}+(1+t)u_1^{t-1}\nabla z_{\alpha}\cdot
\nabla u_1\|\cdot \|u_1\Delta z_{\alpha}+2\nabla z_{\alpha}\cdot
\nabla u_1\|.
\end{aligned}
\end{equation}
For any positive integer $k$, we can find some $\alpha_0$ such that
$$
\sum_{j=\alpha_0+1}^{k}\dfrac{(\lambda_{k+1}-\lambda_j)A_{\alpha_0
j}^2}{\int_{\Omega}|\nabla  z_{\alpha_0}|^2u_1^2}
 = \max_{1\leq \alpha \leq N}{\sum_{j=\alpha+1}^{k}\dfrac{(\lambda_{k+1}-\lambda_j)A_{\alpha j}^2}{\int_{\Omega}|\nabla  z_{\alpha}|^2u_1^2}}.
$$
Hence,  from Lemma 2.1, we obtain
\begin{equation*}
\begin{aligned}
&\dfrac{n(\lambda_{k+1}-\lambda_1)\int_{\Omega}u_1^{t+1}} {1+
\sum\limits_{j=\alpha_0+1}^{k}(\lambda_{k+1}-\lambda_j)\dfrac{A_{\alpha_0
j}^2}
{\int_{\Omega}|\nabla  z_{\alpha_0}|^2u_1^2}}\\
\leq &\sum_{\alpha=1}^N\|u_1^t\Delta z_{\alpha}+(1+t)u_1^{t-1}\nabla
z_{\alpha}\cdot \nabla u_1\|\cdot
\|u_1\Delta z_{\alpha}+2\nabla  z_{\alpha}\cdot \nabla u_1\| \\
\leq &\sqrt{\biggl(n^2\sup_{\Omega}
|H|^2+4\lambda_1\biggl)\biggl(n^2\sup_{\Omega}
|H|^2+\dfrac{(1+t)^2}{2t-1}\lambda_1\biggl)\int_{\Omega}u_1^{2t}},
\end{aligned}
\end{equation*}
that is, we have
\begin{equation}
\begin{aligned}
&\dfrac{n(\lambda_{k+1}-\lambda_1)} {1+
\sum\limits_{j=\alpha_0+1}^{k}(\lambda_{k+1}-\lambda_j)
\dfrac{A_{\alpha_0 j}^2}{\int_{\Omega}|\nabla  z_{\alpha_0}|^2u_1^2}}\\
\leq & ~B(t)\sqrt{\left(n^2\sup_{\Omega}
|H|^2+4\lambda_1\right)\left(n^2\sup_{\Omega}
|H|^2+\dfrac{(1+t)^2}{2t-1}\lambda_1\right)}.
\end{aligned}
\end{equation}

On the other hand, we have
\begin{equation}
\begin{aligned}
&\int_{\Omega}|\nabla  u_1^{t-1}|^2u_1^2
=\dfrac{(t-1)^2}{2t-1}\int_{\Omega}\nabla  u_1\cdot \nabla
u_1^{2t-1} =\dfrac{(t-1)^2}{2t-1}\lambda_1\int_{\Omega}u_1^{2t}.
\end{aligned}
\end{equation}
 Letting
$$
D_j=\int_{\Omega}u_1^{t }u_j,
$$
we know
\begin{equation}
\begin{aligned}
u^t_1=\sum_{j=1}^{\infty}D_{j}u_j, \ \  \  \int_{\Omega}u_1^{2t}=\sum_{j=1}^{\infty}D_{j}^2.
\end{aligned}
\end{equation}
Taking $f=u_1^{t-1}$ in (2.6), we get
\begin{equation*} \aligned \int_\Omega |\nabla
u_{1}^{t-1}|^{2}u_1^{2}&=-2\int_{\Omega}u_1^{t} \nabla  u_1^{t-1}
\cdot \nabla u_1 -\int_{\Omega}u_1^{t+1}\Delta
u_1^{t-1}\\
&=-\sum_{j=1}^{\infty}D_j\biggl(2\int_{\Omega}u_j \nabla  u_1^{t-1}
\cdot \nabla u_1+\int_{\Omega}u_ju_1\Delta
u_1^{t-1}\biggl)\\
&=-\sum_{j=1}^{\infty}D_j\biggl(\int_{\Omega}u_j \Delta
u_1^{t}-\int_{\Omega}u_ju_1^{t-1}\Delta
u_1\biggl)\\
&=-\sum_{j=1}^{\infty}D_j\biggl(\int_{\Omega} u_1^{t} \Delta
u_j-\int_{\Omega}u_ju_1^{t-1}\Delta
u_1\biggl)\\
&=\sum_{j=1}^{\infty}D_j\biggl(\lambda_j\int_{\Omega} u_1^{t}
u_j-\lambda_1\int_{\Omega}u_ju_1^{t}\biggl)\\
&=\sum_{j=2}^{\infty}(\lambda_j-\lambda_1)D_j^{2}.\\
\endaligned
\end{equation*}
 Thus, we infer
\begin{equation}
\begin{aligned}
\sum_{j=2}^{\infty}(\lambda_j-\lambda_1)D_{j}^2=\dfrac{(t-1)^2}{2t-1}\lambda_1\int_{\Omega}u_1^{2t}.
\end{aligned}
\end{equation}
Defining
\begin{equation}
\begin{aligned}
\beta_j=\dfrac{D_{j}}{\sqrt{\dfrac{(t-1)^2}{2t-1}\lambda_1\int_{\Omega}u_1^{2t}}},
\end{aligned}
\end{equation}
we have
\begin{equation}
\begin{aligned}
\sum_{j=2}^{\infty}(\lambda_j-\lambda_1)\beta_{j}^2=1.
\end{aligned}
\end{equation}
For any positive integer $l$,
\begin{equation*}
\begin{aligned}
1=&\sum_{j=2}^{\infty}(\lambda_j-\lambda_1)\beta_{j}^2\\
=&\sum_{j=2}^{l}(\lambda_j-\lambda_1)\beta_{ j}^2+\sum_{j=l+1}^{\infty}(\lambda_j-\lambda_1)\beta_{ j}^2\\
\geq &\sum_{j=2}^{l}(\lambda_j-\lambda_1)\beta_{ j}^2
+(\lambda_{l+1}-\lambda_1)\sum_{j=l+1}^{\infty}\beta_{ j}^2\\
=&\sum_{j=2}^{l}(\lambda_j-\lambda_1)\beta_{ j}^2
+(\lambda_{l+1}-\lambda_1)\sum_{j=2}^{\infty}\beta_{ j}^2-(\lambda_{l+1}-\lambda_1)\sum_{j=2}^{l}\beta_{ j}^2\\
 =&\sum_{j=2}^{l}(\lambda_j-\lambda_{l+1})\beta_{ j}^2
+(\lambda_{l+1}-\lambda_1)\sum_{j=2}^{\infty}\beta_{ j}^2,
\end{aligned}
\end{equation*}
namely,
\begin{equation}
\begin{aligned}
(\lambda_{l+1}-\lambda_1)\sum_{j=2}^{\infty}\beta_{ j}^2\leq
1+\sum_{j=2}^{l}(\lambda_{l+1}-\lambda_j)\beta_{ j}^2.
\end{aligned}
\end{equation}
From (2.29) and (2.30), we infer
\begin{equation}
\begin{aligned}
\dfrac{(t-1)^2}{2t-1}\lambda_1\sum_{j=1}^{\infty}\beta_{ j}^2=1.
\end{aligned}
\end{equation}
Since
\begin{equation*}
\begin{aligned}&\biggl(\int_{\Omega}u_1^{t +1}\biggl)^2
=D_{1}^2=\dfrac{(t-1)^2}{2t-1}\lambda_1\beta_1^2\int_{\Omega}u_1^{2t},
\end{aligned}
\end{equation*}
according to the definition of $B(t)$, we have
\begin{equation}
\begin{aligned}&
B(t)^2=\dfrac{\int_{\Omega}u_1^{2t}}{\biggl(\int_{\Omega}u_1^{t
+1}\biggl)^2} =\dfrac1{\dfrac{(t-1)^2}{2t-1}\lambda_1\beta_1^2}=
\dfrac1{1-\dfrac{(t-1)^2}{2t-1}\lambda_1\sum\limits_{j=2}^{\infty}\beta_{
j}^2}.
\end{aligned}
\end{equation}
From (2.27), (2.33)  and (2.35), we have
\begin{proposition}
Let $M$ be an $n$-dimensional complete  submanifold in  $\mathbf
R^N$, $\Omega \subset M $  a bounded domain with a piecewise smooth
boundary
 $\partial \Omega$. Then,  for any positive integer $k$,
 there exists an integer $\alpha_0$ with $1\leq \alpha_0\leq N$ such that  eigenvalues of
the Dirichlet eigenvalue problem of the Laplacian satisfy,  for any
positive integer $l$  and $t >\frac12$,
\begin{equation*}
\dfrac{n(\lambda_{k+1}-\lambda_1)} {1+
\sum\limits_{j=\alpha_0+1}^{k}(\lambda_{k+1}-\lambda_j)
\dfrac{A_{\alpha_0 j}^2}{\int_{\Omega}|\nabla z_{\alpha_0}|^2u_1^2}}
\leq \sqrt{\dfrac{\biggl(n^2\sup\limits_{\Omega}
|H|^2+4\lambda_1\biggl)\biggl(n^2\sup\limits_{\Omega}
|H|^2+\dfrac{(1+t)^2}{2t-1}\lambda_1\biggl)}{1-\dfrac{(t-1)^2}{2t-1}\dfrac{\lambda_1}{\lambda_{l+1}-\lambda_1}
\biggl(1+\sum\limits_{j=2}^{l}(\lambda_{l+1}-\lambda_j)\beta_{
j}^2\biggl)}}.
\end{equation*}
\end{proposition}

\section{Estimates for eigenvalues on Minimal submanifolds}

In this section, we will deal with eigenvalues of the Laplacian on
bounded domains in complete minimal submanifolds of Euclidean
spaces. Thus, let $\Omega \subset M$ be a bounded domain with
a piecewise smooth boundary
 $\partial \Omega$    in an
 $n$-dimensional complete minimal submanifold $M$ of the Euclidean space $\mathbf R^N$.
 We consider the
 following Dirichlet eigenvalue problem of the Laplacian:
 $$
\left \{ \aligned \Delta u&=-\lambda u \quad \text{in \ \ $\Omega$}, \\
u&=0 \quad \ \ \ \ \text{on $\partial \Omega$}. \endaligned \right.
$$
Since  $M$ is an $n$-dimensional complete  minimal submanifold in
$\mathbf R^N$, we have from Lemma 2.1 and the definition of
$C_{\alpha j}$,
$$
C_{\alpha j}=0
$$
for any $\alpha$ and $j$.
Hence, we have from (2.2),
\begin{equation}
2B_{\alpha j}=(\lambda_1-\lambda_j)A_{\alpha j}.
\end{equation}
For any $\alpha$, we have
\begin{equation*}
\begin{aligned}
0&=-\frac 2{t+1}\int_{\Omega}u_1^{t +1}\Delta z_{\alpha}\\
&= 2\int_{\Omega}u_1^{t}\nabla z_{\alpha}\cdot\nabla u_1\\
&=2\sum_{i=1}^{\infty}D_{i}B_{\alpha
i}=\sum_{i=1}^{\infty}(\lambda_1-\lambda_i)D_{i}A_{\alpha i}.
\end{aligned}
\end{equation*}
Thus, from (2.31) we obtain
\begin{equation}
\sum_{i=2}^{\infty}(\lambda_i-\lambda_1)\beta_{i}A_{\alpha i}=0.
\end{equation}
For any positive integer $j\geq2$, since
\begin{equation*}
\biggl((\lambda_{j}-\lambda_1)\beta_{j}A_{\alpha j}\biggl)^2 \leq
\biggl(\sum_{i=2, i\neq
j}^{\infty}(\lambda_i-\lambda_1)\beta_{i}^2\biggl)\biggl(\sum_{i=2,
i\neq j}^{\infty}(\lambda_i-\lambda_1)A_{\alpha i}^2\biggl),
\end{equation*}
according to (2.8) and (2.32), we derive
\begin{equation*}
(\lambda_{j}-\lambda_1)^2\beta_{j}^2A_{\alpha j}^2 \leq
\biggl(1-(\lambda_{j}-\lambda_1)\beta_{j}^2\biggl)
\biggl(\int_{\Omega}|\nabla  z_{\alpha}|^2u_1^2-
(\lambda_{j}-\lambda_1)A_{\alpha j}^2\biggl).
\end{equation*}
Hence, we have
\begin{equation}
(\lambda_{j}-\lambda_1)\beta_{j}^2+
(\lambda_{j}-\lambda_1)\dfrac{A_{\alpha j}^2}{\int_{\Omega}|\nabla
z_{\alpha}|^2u_1^2}\leq 1.
\end{equation}

From (2.9) and Lemma 2.1, we can get
\begin{equation}
(\lambda_{k+1}-\lambda_{1})\sum_{\alpha=1}^{N}\|z_{\alpha}u_{1}\|^{2}\leq
n+\sum_{\alpha=1}^{N}\sum_{j=\alpha+1}^{k}
(\lambda_{k+1}-\lambda_{j})A_{\alpha j}^{2}.\end{equation}
 On the other hand, from the Stokes theorem,  we have
 $$\aligned
 \int_{\Omega}u_{1}^{t+1}|\nabla z_{\alpha}|^{2}&=\frac{1}{2}\int_{\Omega}u_{1}^{t+1}\Delta z_{\alpha}^{2}
 =-(t+1)\int_{\Omega} z_{\alpha}u_{1}^{t}\nabla
z_{\alpha}\cdot\nabla u_{1}.\\\endaligned$$

 From Lemma 2.1 and the
Cauchy-Schwarz inequality, we get
$$\aligned n\int_{\Omega}u_{1}^{t+1}&=-(t+1)\sum_{\alpha=1}^{N}\int_{\Omega} z_{\alpha}u_{1}^{t}\nabla
z_{\alpha}\cdot\nabla u_{1}\\
&\leq (t+1)\left(\sum_{\alpha=1}^{N}\|
z_{\alpha}u_{1}\|^{2}\right)^{\frac{1}{2}}
\left(\sum_{\alpha=1}^{N}\int_{\Omega} u_1^{2t-2} (\nabla
z_{\alpha}\cdot\nabla u_{1})^{2}\right)^{\frac{1}{2}}\\
&= (t+1)\left(\sum_{\alpha=1}^{N}\|
z_{\alpha}u_{1}\|^{2}\right)^{\frac{1}{2}} \left(\int_{\Omega}
u_1^{2t-2} |\nabla u_{1}|^{2}\right)^{\frac{1}{2}}\\
&= (t+1)\left(\sum_{\alpha=1}^{N}\|
z_{\alpha}u_{1}\|^{2}\right)^{\frac{1}{2}}
\left(\frac{1}{2t-1}\int_{\Omega} \nabla u_{1}^{2t-1}\cdot\nabla
u_1\right)^{\frac{1}{2}}\\
&= (t+1)\left(\sum_{\alpha=1}^{N}\|
z_{\alpha}u_{1}\|^{2}\right)^{\frac{1}{2}}
\left(\frac{\lambda_{1}}{2t-1}\int_{\Omega}
u_{1}^{2t}\right)^{\frac{1}{2}},
\endaligned$$
namely,
\begin{equation}
\sum_{\alpha=1}^{N}\|
z_{\alpha}u_{1}\|^{2}\geq\dfrac{n^{2}\left(\int_{\Omega}u_{1}^{t+1}\right)^{2}}
{\dfrac{(t+1)^{2}}{2t-1}\lambda_{1}\int_{\Omega}u_{1}^{2t}}
=\dfrac{n^{2}}{\dfrac{(t+1)^{2}}{2t-1}\lambda_{1}B(t)^{2}}.
\end{equation}
From (3.3)-(3.5), (2.35) and (2.33), we have
\begin{equation*}\aligned
&\frac{n^{2}(\lambda_{k+1}-\lambda_{1})}{n+\sum\limits_{\alpha=1}^{N}\sum\limits_{j=\alpha+1}^{k}
(\lambda_{k+1}-\lambda_{j})A_{\alpha j}^{2}}\\
 \leq &~
\frac{(t+1)^{2}}{2t-1}\lambda_{1}B(t)^{2}\\
=&~\dfrac{(t+1)^{2}}{2t-1}\dfrac{\lambda_{1}}{1-\dfrac{(t-1)^2}{2t-1}\lambda_1\sum\limits_{j=2}^{\infty}\beta_{
j}^2}\\
\leq
&~\dfrac{(t+1)^{2}}{2t-1}\dfrac{\lambda_{1}}{1-\dfrac{(t-1)^2}{2t-1}\dfrac{\lambda_1}{\lambda_{l+1}-\lambda_{1}}
\biggl(1+\sum\limits_{j=2}^{l}(\lambda_{l+1}-\lambda_{j})\beta_{
j}^2\biggl)}\\
\leq
&~\dfrac{(t+1)^{2}}{2t-1}\dfrac{\lambda_{1}}{1-\dfrac{(t-1)^2}{2t-1}\dfrac{\lambda_1}{\lambda_{l+1}-\lambda_{1}}
\biggl(1+\sum\limits_{j=2}^{l}\dfrac{\lambda_{l+1}-\lambda_{j}}{\lambda_{j}-\lambda_{1}}\biggl[1-(\lambda_{j}-\lambda_{1})
\dfrac{A_{\alpha j}^2}{\int_{\Omega}|\nabla
z_{\alpha}|^2u_1^2}\biggl]\biggl)}.
\endaligned
\end{equation*}

Defining
$$
\sigma_{\alpha l}=\lambda_1+\dfrac{\lambda_{l+1}-\lambda_1}
{1+\sum\limits_{j=2}^{l}\dfrac{\lambda_{l+1}-\lambda_j}{\lambda_{j}-\lambda_1}
\biggl[1-(\lambda_{j}-\lambda_1)\dfrac{A_{\alpha
j}^2}{\int_{\Omega}|\nabla  z_{\alpha}|^2u_1^2}\biggl]}
$$
and taking
$$
t=\dfrac{2\sigma_{\alpha l}}{\sigma_{\alpha l}+\lambda_1},
$$
we obtain the following:
\begin{theorem}
Let $M$ be an $n$-dimensional complete  minimal submanifold in
$\mathbf R^N$, $\Omega \subset M $  a bounded domain with a piecewise
smooth boundary
 $\partial \Omega$. Then,  for positive integers $k, \ l,$ eigenvalues of
the Dirichlet eigenvalue problem of the Laplacian satisfy, for
$1\leq\alpha\leq N,$
\begin{equation}
\begin{aligned}
&\frac{n^{2}(\lambda_{k+1}-\lambda_1)} {n+\sum\limits_{\alpha=1}^{N}
\sum\limits_{j=\alpha+1}^{k}(\lambda_{k+1}-\lambda_j) A_{\alpha
j}^2}
 \leq
3\lambda_1+\dfrac{\lambda_{1}^{2}}{\sigma_{\alpha l}}.
\end{aligned}
\end{equation}
\end{theorem}

\begin{corollary}
Let $M$ be an $n$-dimensional complete  minimal submanifold in
$\mathbf R^N$, $\Omega \subset M $  a bounded domain with a piecewise
smooth boundary
 $\partial \Omega$. Then,  for
the Dirichlet eigenvalue problem of the Laplacian, we have
\begin{equation}
\dfrac{\lambda_{2}} {\lambda_1} \leq \dfrac{n+3+\sqrt{n^{2}+10n+9}}
 {2n}.
\end{equation}
\end{corollary}

\proof\ Taking $k=l=1$ in (3.6), we have
$$n(\lambda_2-\lambda_1)\leq
3\lambda_1+\frac{\lambda_1^{2}}{\lambda_2}.$$ The above inequality
can be written by the following quadratic inequality:
$$
n\left(\frac{\lambda_2}{\lambda_1}\right)^{2}
-(n+3)\frac{\lambda_2}{\lambda_1}-1\leq 0.
$$ Therefore, we can
obtain (3.7).
 \endproof

 \begin{remark}
When $n=2$, the inequality (3.7) becomes the following form:
$$\frac{\lambda_2}{\lambda_1}\leq\frac{5+\sqrt{33}}{4}.$$ Thus, the result of
Brands \cite{B} for a bounded domain in the Euclidean space is also included here.
\end{remark}

 For any positive integer $k$, we can  find  some $\alpha_0$ such that
$$
\sum_{j=\alpha_0+1}^{k}\dfrac{(\lambda_{k+1}-\lambda_j)A_{\alpha_0
j}^2}{\int_{\Omega}|\nabla  z_{\alpha_0}|^2u_1^2}
 = \max_{1\leq \alpha \leq N}{\sum_{j=\alpha+1}^{k}\dfrac{(\lambda_{k+1}-\lambda_j)A_{\alpha j}^2}{\int_{\Omega}|\nabla  z_{\alpha}|^2u_1^2}}.
$$
Then, from Lemma 2.1, we get
$$\aligned &n+\sum_{\alpha=1}^{N}
\sum_{j=\alpha+1}^{k}(\lambda_{k+1}-\lambda_j) A_{\alpha j}^2\\
\leq ~& n+\sum_{j=\alpha_0+1}^{k}(\lambda_{k+1}-\lambda_j)
\frac{A_{\alpha_0 j}^2}{\int_{\Omega}|\nabla
z_{\alpha_0}|^2u_1^2}\sum_{\alpha=1}^{N}\int_{\Omega}|\nabla
z_{\alpha}|^2u_1^2\\
=~&n\biggl(1+\sum_{j=\alpha_0+1}^{k}(\lambda_{k+1}-\lambda_j)
\frac{A_{\alpha_0 j}^2}{\int_{\Omega}|\nabla
z_{\alpha_0}|^2u_1^2}\biggl).
\endaligned$$
Therefore, we have the following

\begin{corollary}
Let $M$ be an $n$-dimensional complete  minimal submanifold in
$\mathbf R^N$, $\Omega \subset M $  a bounded domain with a  piecewise
smooth boundary
 $\partial \Omega$. Then,  for any positive integer $k$,
 there exists an integer $\alpha_0$ with $1\leq \alpha_0\leq N$ such that  eigenvalues of
the Dirichlet eigenvalue problem of the Laplacian satisfy, for any
positive integer $l$,
\begin{equation}
\begin{aligned}
&\dfrac{n(\lambda_{k+1}-\lambda_1)} {1+
\sum\limits_{j=\alpha_0+1}^{k}(\lambda_{k+1}-\lambda_j)
\dfrac{A_{\alpha_0 j}^2}{\int_{\Omega}|\nabla z_{\alpha_0}|^2u_1^2}}
 \leq
3\lambda_1+\dfrac{\lambda_{1}^{2}}{\sigma_{\alpha_{0} l}}.
\end{aligned}
\end{equation}
\end{corollary}

Since  (3.3) holds for any $j$ and any $\alpha$,  from
 Corollary
3.2, we have
\begin{corollary}
Let $M$ be an $n$-dimensional complete  minimal submanifold in
$\mathbf R^N$, $\Omega \subset M $  a bounded domain with a piecewise
smooth boundary
 $\partial \Omega$. Then,  for positive integers $k, \ l$, eigenvalues of
the Dirichlet eigenvalue problem of the Laplacian satisfy, for
$1\leq\alpha\leq N,$
\begin{equation*}
\begin{aligned}
&\dfrac{n(\lambda_{k+1}-\lambda_1)} {1+
\sum\limits_{j=2}^{k}\dfrac{\lambda_{k+1}-\lambda_j}{\lambda_j-\lambda_1}}
\leq \dfrac{n(\lambda_{k+1}-\lambda_1)}
 {1+ \sum\limits_{j=2}^{k}(\lambda_{k+1}-\lambda_j)\biggl(\dfrac{1}{\lambda_j-\lambda_1}-\beta_j^2\biggl)}
\leq 3\lambda_1+\dfrac{\lambda_{1}^{2}}{\sigma_{\alpha l}}.
\end{aligned}
\end{equation*}
\end{corollary}


\providecommand{\bysame}{\leavevmode\hbox
to3em{\hrulefill}\thinspace}

\end{document}